\documentclass[11pt,a4paper]{amsart}
\usepackage{graphics,epic}
\usepackage{amsmath,amssymb, amsthm}
\usepackage[all,2cell]{xy}
\usepackage{xcolor}
\usepackage{textcomp}

\newtheorem{theorem}{Theorem}[section]
\newtheorem*{theorem*}{Theorem}

\newtheorem{conjecture}[theorem]{Conjecture}
\newtheorem*{conjecture*}{Conjecture}
\newtheorem{question}[theorem]{Question}
\newtheorem{example}[theorem]{Example}

\newcommand{\ie}{{\em i.e.}\ }
\newcommand{\confer}{{\em cf.}\ }

\newcommand{\ko}{\: , \;}

\newcommand{\opname}[1]{\operatorname{\mathsf{#1}}}

\renewcommand{\mod}{\opname{mod}\nolimits}

\newcommand{\Mod}{\opname{Mod}\nolimits}

\newcommand{\per}{\opname{per}\nolimits}
\newcommand{\add}{\opname{add}\nolimits}

\newcommand{\thick}{\opname{thick}\nolimits}

\newcommand{\Hom}{\opname{Hom}}

\newcommand{\RHom}{\opname{RHom}}

\newcommand{\End}{\opname{End}}

\newcommand{\ten}{\otimes}
\newcommand{\lten}{\overset{\boldmath{L}}{\ten}}

\newcommand{\ca}{{\mathcal A}}

\newcommand{\cc}{{\mathcal C}}
\newcommand{\cd}{{\mathcal D}}
\newcommand{\ce}{{\mathcal E}}

\newcommand{\cs}{{\mathcal S}}

\renewcommand{\hat}[1]{\widehat{#1}}

\newcommand{\del}{\partial}

\newcommand{\MCM}{\opname{MCM}\nolimits}

\numberwithin{equation}{section}

\begin{document}

\title[2-3-CY]{The interplay between 2-and-3-Calabi--Yau triangulated categories}

\author{Dong Yang}
\address{Dong Yang, Department of Mathematics, Nanjing University, Nanjing 210093, P. R. China}
\email{yangdong@nju.edu.cn}
\thanks{The Chinese version of this survey is published in Sci. Sin. Math.: Yang D. The interplay between 2-Calabi-Yau and 3-Calabi-Yau triangulated categories (in Chinese). Sci Sin Math, 2018, 48: 1Ð10.}
\date{Last modified on \today}

\begin{abstract} 
This short note surveys the constructions of 3-Calabi--Yau triangulated categories with simple-minded collections due to Ginzburg and Kontsevich--Soibelman and the constructions of 2-Calabi--Yau triangulated categories with cluster-tilting objects due to Buan--Marsh--Reineke--Reiten--Todorov and Amiot, and includes a discussion on the normal form of 2-Calabi--Yau triangulated categories with cluster-tilting objects.
\end{abstract}

\maketitle

\tableofcontents

This is a short survey on some `recent' progresses on the interplay between 2-Calabi--Yau and 3-Calabi--Yau triangulated categories. It is based on my talks given in Stuttgart in November 2011 and in Nagoya in May 2013. 

2-Calabi--Yau triangulated categories with cluster-tilting objects play a central role in the additive categorification of Fomin--Zelevinsky's cluster algebras, see for example \cite{BuanMarshReinekeReitenTodorov06,BuanMarshReiten08,CalderoKeller08,Palu08,Plamondon11}. A prototypical example of such categories is the cluster category of an acyclic quiver, introduced by Buan--Marsh--Reineke--Reiten--Todorov in \cite{BuanMarshReinekeReitenTodorov06}.

3-Calabi--Yau triangulated categories with simple-minded collections play an important role in algebraic geometry and mathematical physics \cite{Bridgeland07a,KontsevichSoibelman08}. There are two constructions of such categories: for a quiver with potential, Ginzburg constructs a 3-Calabi--Yau dg algebra $\Gamma$ \cite{Ginzburg06}, whose finite-dimensional derived category $\cd_{fd}(\Gamma)$ is a 3-Calabi--Yau triangulated category with simple-minded collections; Kontsevich--Soibelman constructs a 3-Calabi--Yau $A_\infty$-algebra $\ca$ \cite{KontsevichSoibelman08}, whose perfect derived category $\per(\ca)$ is a 3-Calabi--Yau triangulated category with simple-minded collections. These two constructions are Koszul-dual to each other.

A direct connection between these two types of categories was observed in \cite{KellerReiten08,Tabuada07,Amiot09,Keller11a}. For a quiver with potential, denote by $\Gamma$ the associated Ginzburg dg algebra. It is shown in \cite{Amiot09,Keller11a} that the generalised cluster category
\[
\cc:=\per(\Gamma)/\cd_{fd}(\Gamma)
\]
is a 2-Calabi--Yau triangulated category with cluster-tilting objects. It is conjectured that in characteristic zero all 2-Calabi--Yau triangulated category with cluster-tilting objects are of this form. This conjecture holds true for all known examples \cite{KellerReiten08,Amiot09,AmiotReitenTodorov11,AmiotIyamaReitenTodorov12,AmiotIyamaReiten15,deVolcseyVandenBergh16}.

This short note surveys the above constructions due to Buan--Marsh--Reineke--Reiten--Todorov \cite{BuanMarshReinekeReitenTodorov06}, Ginzburg \cite{Ginzburg06}, Kontsevich--Soibelman \cite{KontsevichSoibelman08}, Amiot \cite{Amiot09} and Keller \cite{Keller11a}, and discusses various methods to attack the conjecture that all 2-Calabi--Yau triangulated categories with cluster-tilting objects are generalised cluster categories of quivers with potential.

\smallskip
Throughout, $k$ is an algebraically closed field and $D=\Hom_k(?,k)$ denotes the $k$-dual. All categories are $k$-categories and all algebras are $k$-algebras.

\section{Silting objects, simple-minded collections and cluster-tilting objects}

Let $\cc$ be a triangulated category with suspension functor $\Sigma$. For a set $\cs$ of objects of $\cc$, let $\thick(\cs)$ denote the smallest thick subcategory of $\cc$ containing $\cs$.

An object $M$ of $\cc$ is a \emph{silting object} (\cite{KellerVossieck88,AssemSoutoTrepode08,AiharaIyama12}) of $\cc$ if 
\begin{itemize}
\item[-]
$\Hom_{\cc}(M,\Sigma^p M)=0$ for any $p>0$, 
\item[-]
$\cc=\thick(M)$.
\end{itemize} 
A collection $\{X_1,\ldots,X_n\}$ of objects of $\cc$ is \emph{simple-minded} (\cite{KoenigYang14,KellerNicolas13}) if
\begin{itemize}
\item[-] $\Hom_{\cc}(X_i,\Sigma^p X_j)=0,~~\forall~p<0$,
\item[-] $\Hom_{\cc}(X_i,X_j)=\begin{cases} k & \text{if\ }i=j,\\
                                           0 & \text{otherwise},
                                           \end{cases}$
\item[-] $\cc=\thick(X_1,\ldots,X_n)$.
\end{itemize}

Fix $d\in\mathbb{Z}$. Assume further that $\cc$ Hom-finite and Krull--Schmidt. 
We say that $\cc$ is \emph{$d$-Calabi--Yau} if 
 there is a bifunctorial isomorphism\[
D\Hom_{\cc}(X,Y)\stackrel{\cong}{\longrightarrow} \Hom_{\cc}(Y,\Sigma^d X)
\]
for any $X,Y\in\cc$ (\cite{Kontsevich98,Keller10}). Note that in \cite{Keller10} such categories are called weakly $d$-Calabi--Yau triangulated categories.
 
Assume that $\cc$ is $d$-Calabi--Yau. An object $T$ of $\cc$ is \emph{$d$-cluster-tilting} (\cite{KellerReiten07,Thomas07,Zhu08,IyamaYoshino08}) if
\begin{itemize}
\item[-] $\Hom_\cc(T,\Sigma^p X)=0$ for $0<p<d$ $\Longleftrightarrow$ $X$ belongs to $\add(T)$.
\end{itemize}
$2$-cluster-tilting objects are usually called cluster-tilting objects.

In this note, we are interested in $3$-Calabi--Yau triangulated categories with simple-minded collections and $2$-Calabi--Yau triangulated categories with cluster-tilting objects.

\section{3-CalabiÐYau triangulated categories}
The algebraic theory of 3-CY's has two sides: the Ginzburg side and the Kontsevich--Soibelman side. Both sides are associated to quivers with potential.

\subsection{The Ginzburg side} On the Ginzburg side there is a dg algebra. 

For a dg algebra $A$, let $\cd(A)$ be the derived category of (right) dg $A$-modules, $\per(A)=\thick(A_A)$ be the perfect derived category, and $\cd_{fd}(A)$ be the finite-dimensional derived category (\ie $M\in\cd_{fd}(A)$ if and only if $H^*(M)$ is finite-dimensional). See \cite{Keller94,Keller06d}.

To a quiver with potential $(Q,W)$, Ginzburg associates a dg algebra $\widehat{\Gamma}(Q,W)$, which we call the \emph{complete Ginzburg dg algebra} of $(Q, W)$, see \cite{Ginzburg06,KellerYang11}. Precisely, $\hat{\Gamma}(Q,W)$ is
constructed as follows: Let $\tilde{Q}$ be the
graded quiver with the same vertices as $Q$ and whose arrows are
\begin{itemize}
\item the arrows of $Q$ (they all have degree~$0$),
\item an arrow $a^*: j \to i$ of degree $-1$ for each arrow $a:i\to j$ of $Q$,
\item a loop $t_i : i \to i$ of degree $-2$ for each vertex $i$
of $Q$.
\end{itemize}
The underlying graded algebra of $\hat{\Gamma}(Q,W)$ is the
completion of the graded path algebra $k\tilde{Q}$ in the category
of graded vector spaces with respect to the ideal generated by the
arrows of $\tilde{Q}$. Thus, the $n$-th component of
$\hat{\Gamma}(Q,W)$ consists of elements of the form
$\sum_{p}\lambda_p p$, where $p$ runs over all paths of degree $n$. The completion endows $\hat{\Gamma}(Q,W)$ with a \emph{pseudo-compact} topology, see \cite{Keller11a,VandenBergh15}. 
The differential of $\hat{\Gamma}(Q,W)$ is the unique continuous
linear endomorphism homogeneous of degree~$1$ which satisfies the graded
Leibniz rule
\[
d(uv)= (du) v + (-1)^p u dv \ko
\]
for all homogeneous $u$ of degree $p$ and all $v$, and takes the
following values on the arrows of $\tilde{Q}$:
\begin{itemize}
\item[-] $d(a)=0$ for each arrow $a$ of $Q$,
\item[-] $d(a^*) = \del_a W$ for each arrow $a$ of $Q$, where $\del_a$ is the cyclic derivative associated to $a$ (see \cite{DerksenWeymanZelevinsky08}; roughly, it removes $a$ from $W$),
\item[-] $d(t_i) = e_i (\sum_{a\in Q_1} [a,a^*]) e_i$ for each vertex $i$ of $Q$, where
$e_i$ is the trivial path at $i$.
\end{itemize}

The \emph{complete Jacobian algebra} of $(Q,W)$ is the $0$-th cohomology of $\hat{\Gamma}(Q,W)$:
\[
\hat{J}(Q,W)=\hat{kQ}/\overline{(\del_a W:a\in Q_1)},
\]
where $\overline{I}$ is the closure of $I$. When $W$ is a finite sum, one can define the non-complete Ginzburg dg algebra and Jacobian algebra.

\begin{theorem} [{\cite[Theorem A.17]{Keller11a},~\cite{Keller11}}] \label{t:Keller-cy} Let $(Q,W)$ be a quiver with potential and let $\Gamma=\hat{\Gamma}(Q,W)$ be the complete Ginzburg dg algebra.

\begin{itemize}
\item[(a)] $\Gamma$ is topologically homologically smooth, \ie $\Gamma$ as a dg $\Gamma$-bimodule belongs to $\per(\Gamma^{op}\hat{\ten}\Gamma)$. In particular, the perfect derived category $\per(\Gamma)$ contains the finite-dimensional derived category $\cd_{fd}(\Gamma)$.
\item[(b)] $\Gamma$ is bimodule 3-Calabi--Yau, \ie there is an isomorphism $\RHom_{\Gamma^{op}\hat{\ten}\Gamma}(\Gamma,\Gamma^{op}\hat{\ten}\Gamma)\cong \Sigma^{-3}\Gamma$ in $\cd(\Gamma^{op}\hat{\ten}\Gamma)$.
\end{itemize}
\end{theorem}

It follows that the triangulated category $\cd_{fd}(\Gamma)$ is 3-Calabi--Yau (\cite[Section A.15]{Keller11a}). Moreover, one checks that the simple $\hat{J}(Q,W)$-modules form a simple-minded collection of $\cd_{fd}(\Gamma)$.

Notice that $\Gamma=\hat{\Gamma}(Q,W)$ is concentrated in non-positive degrees. It follows that $\Gamma$ is a silting object of $\per(\Gamma)$. We have the converse of Theorem \ref{t:Keller-cy}, due to Van den Bergh, as a special case of \cite[Theorems A and B]{VandenBergh15}.

\begin{theorem}[\cite{VandenBergh15}] Assume that $k$ is of characteristic $0$.
Let $A$ be a topologically homologically smooth bimodule $3$-Calabi--Yau pseudo-compact dg algebra such that $A$ is a silting object in $\per(A)$. Then $A$ is quasi-isomorphic to the complete Ginzburg dg algebra of some quiver with potential.
\end{theorem}

\subsection{The Kontsevich--Soibelman side} On the Kontsevich--Soibelman side, there is an $A_\infty$-algebra. An $A_\infty$-algebra $A$ is a graded vector space endowed with a family of maps $m_n: A^{\ten n}\rightarrow A$ homogeneous of degree $2-n$ satisfying certain conditions, see for example \cite{Keller06c}. For an $A_\infty$-algebra $A$, one can define the derived category $\cd(A)$, the perfect derived category $\per(A)$ and the finite-dimensional derived category $\cd_{fd}(A)$ as well.

A cyclic structure of degree $d$ on an $A_\infty$-algebra $A$ is a supersymmetric non-degenerate bilinear form of degree $d$
\[
(-,-): A\times A\rightarrow k[-d]
\]
such that 
\[
(m_n(a_1,\ldots,a_n), a_{n+1}) = (-1)^n(-1)^{|a_1|(|a_2|+\ldots+|a_{n+1}|)}(m_n(a_2,\ldots, a_{n+1}), a_1).
\]

To a quiver with potential Kontsevich--Soibelman associates an $A_\infty$-algebra $\ca(Q,W)$, which we call the \emph{Kontsevich--Soibelman $A_\infty$-algebra}. Precisely, as a graded vector space $\ca(Q,W)$ has a basis concentrated in degrees $0,1,2,3$:
\begin{itemize}
\item the trivial path $e_i$ of $Q$ in degree $0$ for each vertex $i$ of $Q$,
\item an element $a^*$ in degree $1$ for each arrow $a$ of $Q$, 
\item the arrows of $Q$, in degree $2$,
\item an element $e_i^*$ in degree $3$ for each vertex $i$ of $Q$.
\end{itemize}
Write $W=\sum_c \lambda_c c$, where the sum is over all cycles $c$ of $Q$ of length $\geq 3$. 
The maps $m_n$ are given by
\begin{itemize}
\item[--] $m_1=0$,
\item[--] $m_2(e_j\ten a)=a=m_2(a\ten e_i)$ and $m_2(e_i\ten a^*)=a^*=m_2(a^*\ten e_j)$ if $a:i\to j$ is an arrow of $Q$,
\item[--] $m_2(e_i\ten e_j^*)=0=m_2(e_j^*\ten e_i)$ and $m_2(e_i\ten e_i^*)=e_i^*=m_2(e_i^*\ten e_i)$ if $i$ and $j$ are different vertices of $Q$,
\item[--] $m_2(a\ten a^*)=t_j^*$ and $m_2(a^*\ten a)=t_i^*$ if $a:i\to j$ is an arrow of $Q$,
\item[--] $m_n(\cdots\ten e_i\ten \cdots)=0$ for any vertex $i$ of $Q$, if $n\geq 3$,
\item[--] $m_n(a_1^*\ten\cdots\ten a_n^*)=-\sum_{a}\lambda_{a_n\cdots a_1 a}a$, where the sum is over all arrows of $Q$.
\end{itemize}

\begin{theorem} [\cite{KontsevichSoibelman08}] \label{t:KS-cyclic}
Assume that $k$ is of charateristic $0$.
\begin{itemize}
\item[(a)] Let $(Q,W)$ be a quiver with potential. Then $\ca(Q,W)$ has a natural cyclic structure of degree $3$.
\item[(b)] Let $A$ be an $A_\infty$-algebra with a cyclic structure of degree $3$ such that the indecomposable direct summands of $A_A$ form a simple-minded collection of $\per(A)$. Then there is a quiver with potential $(Q, W )$ such that $A$ is $A_\infty$-isomorphic to $\ca(Q,W)$.
\end{itemize}
\end{theorem}

\subsection{The connection} The Ginzburg side and the Kontsevich--Soibelman side are related by Koszul duality. Precisely, the complete Ginzburg dg algebra $\Gamma=\widehat{\Gamma}(Q,W)$ of a quiver with potential $(Q,W)$ can be obtained from the Kontsevich--Soibelman $A_\infty$-algebra $\ca=\ca(Q,W)$ by taking the dual bar construction (\cite{Lefevre03,LuPalmieriWuZhang08}):
\[
{\Gamma}=\Hom_K(T_K(\ca^{\geq 1}[1]),K),
\]
where $K$ is the semi-simple algebra $kQ_0$.
Conversely, the $A_\infty$-Koszul dual of $\Gamma$ is the Kontsevich--Soibelman $A_\infty$-algebra $\ca$. This leads to the following triangle equivalences
\[
\xymatrix{
\per(\Gamma)\ar[rr]^{\simeq} && \cd_{fd}(\ca)\\
\cd_{fd}(\Gamma)\ar[rr]^{\simeq}\ar@{^{(}->}[u] && \per(\ca). \ar@{^{(}->}[u]
}
\]

\section{2-Calabi--Yau triangulated categories}
\label{s:2-cy-cat}
A typical 2-CY triangulated category is the cluster category associated to an acyclic quiver.

Let $Q$ be an acyclic quiver, \ie a finite quiver without oriented cycles. Let $\mod kQ$ denote the category of finite-dimensional modules over the path algebra $kQ$ and let $\cd^b(\mod kQ)$ denote the corresponding bounded derived category. The derived Nakayama functor $\nu=?\lten_{kQ} D(kQ)$ is a Serre functor of $\cd^b(\mod kQ)$  (\cite[Theorem 4.6]{Happel87}). Define the \emph{cluster category} $\cc_Q$ (\cite{BuanMarshReinekeReitenTodorov06}) to be the orbit category
\[
\xymatrix{
\cc_Q:=\cd^b(\mod kQ)/\nu\Sigma^{-2}.
}
\]

\begin{theorem}
\begin{itemize}
\item[(a)] \emph{(\cite[Theorem 1 and Corollary 1]{Keller05})} $\cc_Q$ is a 2-CY triangulated category.
\item[(b)] \emph{(\cite[Theorem 3.3(b)]{BuanMarshReinekeReitenTodorov06})} The image of the $kQ$ in $\cc_Q$ is a cluster-tilting object.
\end{itemize}
\end{theorem}

There are many 2-CY triangulated categories with cluster-tilting objects arising from
\begin{itemize}
\item[(1)] module categories of preprojective algebras of acyclic quivers (\cite{GeissLeclercSchroer08a,BuanIyamaReitenScott09}),
\end{itemize}
and from 
\begin{itemize}
\item[(2)]
categories of maximal Cohen--Macaulay modules of singularities (\cite{BurbanIyamaKellerReiten08,Iyama07}).
\end{itemize}

There is the following  remarkable recognition theorem due to Keller and Reiten.

\begin{theorem}[{\cite[Theorem 2.1]{KellerReiten08}}]\label{t:KR-recognition}
Let $\cc$ be a 2-Calabi--Yau algebraic triangulated category with a cluster-tilting object whose endomorphism algebra is the path algebra $kQ$ of an acyclic quiver $Q$. Then there is a triangle equivalence
\[
\cc_Q\stackrel{\simeq}{\longrightarrow} \cc.
\]
\end{theorem}

\section{From 3-CYÕs to 2-CYÕs}
Let $(Q,W)$ be a quiver with potential. The \emph{generalised cluster category} (\cite{Amiot09}) is defined as 
\[
\cc_{(Q,W)} := \per(\widehat{\Gamma}(Q,W))/\cd_{fd}(\widehat{\Gamma}(Q,W)).\]

The relation among the complete Ginzburg dg algebra $\widehat{\Gamma}(Q, W )$, the Kontsevich--Soibelman $A_\infty$-
algebra $\ca(Q,W)$ and the generalised cluster category $\cc_{(Q,W)}$ is encoded in the following recollement 

\[\xymatrix{\widetilde{C}_{(Q,W)}\ar[rr]&&\cd\widehat{\Gamma}{(Q,W)}\ar[rr]\ar@/^20pt/[ll]\ar@/_20pt/[ll]&&\cd\ca{(Q,W)}\ar@/^20pt/[ll]\ar@/_20pt/[ll]},
\vspace{10pt}
\]
where $\tilde{\cc}_{(Q,W )}$ is an unbounded version of $\cc_{(Q,W )}$ in the sense that $\tilde{\cc}_{(Q,W )}$ has infinite direct
sums and $\cc_{(Q,W)}$ consists of compact objects in $\tilde{\cc}_{(Q,W)}$ (\cite[Corollary 3]{Yang12a}).

\begin{theorem} If $\widehat{J}(Q,W)$ is finite-dimensional, then $\cc_{(Q,W )}$ is Hom-finite and 2-CY as a triangulated category. Moreover, the image of the object $\widehat{\Gamma}(Q,W)$ in $\cc_{(Q,W)}$ is a cluster-tilting object  whose endomorphism algebra is $\widehat{J}(Q,W)$.
\end{theorem}

This result is a special case of the following more general theorem, which is a `topological' version of  \cite[Theorems 2.1 and 3.5]{Amiot09}. It is generalised in \cite{Guolingyan11a, IyamaYang18}.

\begin{theorem}[{\cite[Theorem A.21]{Keller11a}}]\label{thm:Amiot}
Let $A$ be a dg algebra satisfying the conditions
 \begin{itemize}
 \item[-] $H^i(A)=0$ for $i>0$,
 \item[-] $H^0(A)$ is finite-dimensional,
 \item[-] $A$ is topologically homologically smooth,
 \item[-] $A$ is bimodule 3-Calabi--Yau.
 \end{itemize}
 Then $\per(A)/\cd_{fd}(A)$ is Hom-finite and 2-CY as a triangulated category. Moreover, the image of $A$ in $\per(A)/\cd_{fd}(A)$ is a cluster-tilting object whose endomorphism algebra is $H^0(A)$.
\end{theorem}

\section{From 2-CYÕs to 3-CYÕs}
Motivated by the Keller--Reiten recognition Theorem~\ref{t:KR-recognition}, we propose the following conjecture (\confer \cite[Summary of results, Part 2, Perspectives]{Amiot08}).

\begin{conjecture}\label{conj:morita-for-2-cy} Assume that $k$ is of characteristic $0$. 
Let $\cc$ be a 2-Calabi--Yau algebraic triangulated category with a cluster-tilting object $T$. Then there is a quiver with potential $(Q,W)$ together with a triangle equivalence
\[
\xymatrix@R=0.5pc{\cc_{(Q,W)}\ar[rr]^\simeq && \cc.\\
\widehat{\Gamma}(Q,W)\ar@{|->}[rr] && T
}
\]
\end{conjecture}

If this conjecture holds, then the answer to the following question proposed in Amiot's ICRA XIV notes is positive.

\begin{question}[{\cite[Question 2.20]{Amiot11}}]
Assume that $k$ is of characteristic $0$. 
Let $\cc$ be a 2-Calabi--Yau algebraic triangulated category with a cluster-tilting object $T$. Then the endomorphism algebra of $T$ is the complete Jacobian algebra of some quiver with potential.
\end{question}

Note that  we have replaced `Jacobian algebra' by `complete Jacobian algebra'. The original question has a negative answer, as there are quivers with potentials whose complete Jacobian algebras are not non-complete Jacobian algebras of any quiver with potential, see \cite[Example 4.3]{Plamondon13} for an example. 

Moreover, we have put an extra assumption on the characteristic of the field, 
as when $k$ is of positive characteristic, the answer to the question is negative. For example, if $k$ is of characteristic $p>0$, then $k[x]/(x^{p-1})$ is not a Jacobian algebra. However, take $\Gamma$ as the dg algebra whose underlying graded algebra is $k\langle\langle x,x^*,t\rangle\rangle$ with $\deg(x)=0$, $\deg(x^*)=-1$ and
$\deg(t)=-2$, and whose differential $d$
is defined by
\[
d(x)=0,~~~ d(x^*)= x^{p-1},~~~ d(t)= xx^*-x^*x.
\]
It is straightforward to check that $\Gamma$ satisfies the assumptions of Theorem~\ref{thm:Amiot}, so $k[x]/(x^{p-1})=H^0\Gamma$ is the endomorphism of a cluster-tilting object in a 2-Calabi--Yau algebraic triangulated category. More examples can be found in  \cite{Ladkani14}.

One approach to attack Conjecture~\ref{conj:morita-for-2-cy} is developed by Amiot \cite{Amiot09}. For many 2-Calabi--Yau triangulated categories $\cc$ arising from (1) and (2) in Section~\ref{s:2-cy-cat}, there is a finite-dimensional algebra $A$ of global dimension $2$ with certain conditions such that $\cc$ is triangle equivalent to the cluster category of $A$. On the other hand, the \emph{3-Calabi--Yau completion} $\Pi_3(A)$ of $A$ in the sense of Keller~\cite{Keller11} satisfies the four conditions in Theorem~\ref{thm:Amiot} (with `topologically homologically smooth' replaced by `homologically smooth'), and the cluster category of $A$ is triangle equivalent to $\per(\Pi_3(A))/\cd_{fd}(\Pi_3(A))$ (\cite[Theorem 6.12(a)]{Keller11}. By \cite[Theorem 6.10]{Keller11}, there is a quiver with potential $(Q,W)$ such that $\Pi_3(A)$ is quasi-isomorphic to the $\Gamma(Q,W)$. Therefore $\cc$ is triangle equivalent to $\per(\Gamma(Q,W))/\cd_{fd}(\Gamma(Q,W))$. 
Amiot, Iyama, Reiten, Todorov \cite{Amiot09,AmiotReitenTodorov11,AmiotIyamaReitenTodorov12,AmiotIyamaReiten15} have intensively worked on this topic.  

\smallskip
Another promising approach is the from-2CY-to-3CY approach, which we explain in more details below.
Let $\ce$ be a stably 2-CY Frobenius category and $M\in\ce$ be a basic cluster-tilting object. Let $\cc$ be the stable category of $\ce$, $A = \End_\ce(M)$ and $e \in A$ be the idempotent corresponding to the maximal projective-injective direct summand of M. Then $A/AeA = \End_{\cc}(M)$.
The following is a consequence of \cite[Proposition 4 (c)]{KellerReiten07}. 

\begin{theorem} \label{thm:cy-localisation}
Let $\cd_{fd,A/AeA}(\Mod A)$ be the subcategory of $\cd(\Mod A)$ whose total cohomology is finite-dimensional and  supported on $A/AeA$. Then there is a bifunctorial isomorphism
\[
D\Hom(X,Y)\stackrel{\cong}{\longrightarrow} \Hom(Y,\Sigma^{d+1} X),
\]
where $X\in \mod A/AeA$ and $Y\in \cd_{fd,A/AeA}(\Mod A)$.
\end{theorem}

The collection of simple $A$-modules  supported on $A/AeA$ is a simple-minded collection in $\cd_{fd,A/AeA}(\Mod A)$. By Lef\`evre's $A_\infty$-version of Morita's theorem for triangulated categories, there is an $A_\infty$-algebra $E$ together with a triangle equivalence $\per(E)\rightarrow \cd_{fd,A/AeA}(\Mod A)$ which takes the direct summands of $E$ to the simple-minded collection above. It follows from Theorem~\ref{thm:cy-localisation} that on $E$ there is a supersymmetric non-degenerated bilinear form. Very sadly, however, we cannot apply Theorem~\ref{t:KS-cyclic}, because it is not clear whether the bilinear form defined a cyclic structure on $E$. 


\smallskip
For $\ce = \MCM(R)$ the category of maximal Cohen--Macauley modules over a complete Gorenstein local algebra $R$ with isolated singularity,  de V\"olcsey and Van den Bergh has another approach in \cite{deVolcseyVandenBergh16}.

\begin{theorem}[{\cite[Theorem 1.1]{deVolcseyVandenBergh16}}] Assume as above with $\ce = \MCM(R)$ for a complete Gorenstein local algebra $R$ with isolated singularity. Assume that $\tilde{A}\rightarrow A$ is a cofibrant minimal model of $A$ (with the idempotent $e$ lifted) and let $B = \tilde{A}/\overline{\tilde{A}e\tilde{A}}$. Then there is a triangle equivalence
$\per(B)/\cd_{fd}(B) \stackrel{\simeq}{\rightarrow} \cc$.
\end{theorem}

It is not hard to check that $\cd_{fd}(B)\simeq \cd_{fd,A/AeA}(\Mod A)$ (\cite[Corollary 2.12(b)]{KalckYang16}). 
If $\tilde{A}$ is the complete Ginzburg dg algebra of some quiver with potential, so is $B$. 
In particular, we obtain a triangle equivalence
\[\cc_{(Q,W)}\cong\cc.\]
Martin Kalck and myself independently started this approach in \cite{KalckYang16,KalckYang14}.

\begin{example}{\rm 
As an example, let us consider the 3-dimensional McKay correspondence. Let $G \subset SL_3(k)$ be a finite subgroup. Let $S = k{\mbox{\textlbrackdbl}} x,y,z\mbox{\textrbrackdbl}$ with the induced $G$-action and let $R = S^G$. Assume that $R$ has isolated singularity. Then by a theorem of Iyama \cite[Theorem 2.5]{Iyama07a}, as an $R$-module $S$ is a 2-cluster-tilting object in $\MCM(R)$. Its endomorphism algebra $A = \End_R(S)$ is isomorphic to $\widehat{J}(Q',W')$, where $Q'$ is the McKay quiver and $W'$ is a ÕcanonicalÕ potential. See \cite[Section 4.4]{Ginzburg06} for the above statement and the construction of $W'$. Moreover, in this case the canonical homomorphism $\widehat{\Gamma}(Q', W')\rightarrow \widehat{J}(Q', W')$ is a quasi-isomorphism, so it is a cofibrant minimal model of $\widehat{J}(Q',W')$. Let $e$ be the identity of $R$, considered as an element of $A$. Let $0$ be the vertex of $Q'$ corresponding to the trivial $G$-module. Let $(Q, W)$ be the quiver with potential obtained from $(Q', W')$ by deleting the vertex $0$. Then $B = \widehat{\Gamma}(Q', W')/\overline{\widehat{\Gamma}(Q', W')e\widehat{\Gamma}(Q', W')} = \widehat{\Gamma}(Q, W)$. As a consequence, the stable category $\underline{\MCM}(R)$ is triangle equivalent to the generalised cluster category $\cc_{(Q,W)}$.
}
\end{example}

\def\cprime{$'$}
\providecommand{\bysame}{\leavevmode\hbox to3em{\hrulefill}\thinspace}
\providecommand{\MR}{\relax\ifhmode\unskip\space\fi MR }
\providecommand{\MRhref}[2]{%
  \href{http://www.ams.org/mathscinet-getitem?mr=#1}{#2}
}
\providecommand{\href}[2]{#2}

\end{document}